\documentclass[a4paper,oneside]{amsart}
\usepackage{geometry}
\usepackage{amsmath}
\usepackage{flexisym}
\usepackage{mathrsfs} 
\usepackage{amsfonts} 
\usepackage[initials, nobysame,shortalphabetic]{amsrefs}
\renewcommand{\parenthesize}[1]{} 
\usepackage{comment}
\usepackage{amssymb}
\usepackage{amsmath}
\usepackage{graphicx} 
\usepackage{hyperref}
\usepackage{todonotes}
\hypersetup{
	colorlinks,
	linkcolor={red!40!black},
	citecolor={red!50!black},
	urlcolor={red!50!black}
}
\usepackage[nameinlink]{cleveref}
\newtheorem{theorem}{Theorem}[section]
\newtheorem{lemma}[theorem]{Lemma}

\newtheorem{proposition}[theorem]{Proposition}
\newtheorem{corollary}[theorem]{Corollary}
\newtheorem{definition}[theorem]{Definition}

\theoremstyle{definition}
\newtheorem{remark}[theorem]{Remark}
\newcommand{\N}{\mathbb{N}}

\newcommand{\Z}{\mathbb{Z}}

\newcommand{\spa}{\operatorname{spa}}
\title{Averages of sub-additive functions through non-F\o lner sequences}
\author{Nicanor Carrasco-Vargas}
\thanks{IRIF, Université Paris Cité, Paris, France.}
\begin{document}
	\maketitle
	\begin{abstract}
		We study the quantity $\sup_{\phi\in\Phi}\limsup_{n\to\infty}\phi(F_n)/|F_n|$, where $\mathcal F= (F_n)_{n\in\N}$ is a sequence of finite subsets of a countable amenable group, and $\Phi$ is a collection of subadditive functions.  Under mild conditions on $\Phi$ we prove that the aforementioned quantity equals the supremum over $\phi\in\Phi$ of the limit of $\phi(F)/|F|$ as $F$ becomes more and more invariant, multiplied by  a ``F\o lnerness'' constant $C(\mathcal F)\in[1,\infty]$. We have $C(\mathcal F)=1$ if and only if $\mathcal F$ is F\o lner. A corollary of this result is that every choice of a non-F\o lner sequence overestimates topological and measure-theoretic entropy by the multiplicative factor $C(\mathcal F)$ (in the positive entropy case).  
		
		We obtain a generalization of classical results of Goodman, Krug, and Newton, on the so-called proportionality principle for sequence entropy, and their extension by Gao, Huang, and Zhu to countable amenable groups. 
	\end{abstract}

	\section{Introduction}
    
    Sequence entropy is an invariant introduced by  
	Kushnirenko \cite{kushnirenko_metric_1967} with the purpose of distinguishing pairs of zero-entropy dynamical systems. Given a sequence $A=(t_i)_{i\in \N}$ of natural numbers, one replaces the iterates $T^1,T^2,\dots$ of a transformation $T$ by $T^{t_1}, T^{t_2},\dots$ in the classical definitions of entropy. The invariant obtained is denoted $h_{\mu}^A(T)$ for a measure-preserving system $(X,\mathcal X,\mu,T)$, and $h_{top}^A(T)$ for a topological dynamical system $(X,T)$. We refer the reader to the survey \cite{kanigowski_survey_2020} and references for historical context, and to \cite{takeda2026sequenceentropyranksystems,zbMATH08015074,gomez2026sequence,zbMATH07802375,zbMATH08095891,zbMATH07779083,zbMATH07862705,zbMATH07854720} for recent developments.

	If one wants to use sequence entropy to distinguish zero-entropy systems, then the choice of the sequence $A$ is crucial. Early investigations showed that if the sequence $A$ is not sufficiently sparse (for instance, a  syndetic sequence) then the invariant obtained is just classical entropy scaled by a multiplicative factor which only depends on the sequence, and in particular is also zero for zero entropy systems. It is also known that sequence entropy does not provide new information for positive entropy systems, in the sense that if we choose a sequence which does not fall into the aforementioned class, then the corresponding sequence entropy is infinite for every positive entropy system.  
	
	The \textit{proportionality principle}  mentioned in the previous paragraph was first discovered by Krug and Newton \cite{krug_sequence_1972} for measure-theoretic sequence entropy, and subsequently generalized by Goodman \cite{goodman_topological_1974} to  topological sequence entropy. A generalization of the same results for actions of countable amenable groups was recently proved by Gao, Huang and Zhu \cite{gao_sequence_2020}, for sequences satisfying a technical condition (see also  \cite{eberlein_topological_1976,zbMATH03820090,zbMATH03586233,huang_combinatorial_2009}). 
	
	Our main result removes this technical condition, and generalizes the previous results to countable amenable groups. Instead of working separately with topological and measure-theoretic sequence entropy, we interpret this proportionality principle as a statement about subadditive functions.

    Let $G$ be a countable amenable group and let $\mathscr F(G)$ be the set of finite nonempty subsets of $G$. A generalization of Fekete's Lemma, often called  Ornstein-Weiss Lemma, states that when $\phi\colon \mathscr{F}(G)\to[0,\infty)$ is subadditive and translation-invariant then 
	\[F\mapsto \frac{1}{|F|}\phi(F)\]
	converges to a limit as $F$ becomes more and more invariant. The next result describes what happens when we compute the limit along a sequence $(F_n)_{n\in\N}$ which is not necessarily F\o lner, and we take the supremum over a family of functions $\phi$ with natural closure properties. 
	\begin{theorem}[Special case of \Cref{thm:main}]\label{thm:main-intro}
		Let $\Phi$ be a family of monotone subadditive and translation-invariant functions from $\mathscr{F}(G)$ to $[0,\infty)$. Suppose that $\Phi$ has the following closure property: for every $\phi\in\Phi$ and $E\Subset G$, the function $F\to\phi(EF)$
		belongs to $\Phi$. For each $\phi\in \Phi$, let $\mathfrak{h}_\phi$ be the limit of $\phi(F)/|F|$ as $F$ becomes more and more invariant, and  define 
		\[\mathfrak{h}=\sup_{\phi\in\Phi}\mathfrak{h}_{\phi}.\]
		Assume further that $\Phi$ satisfies  $\mathfrak{h}=\sup_{\phi\in\Phi}\inf_{F\Subset G}\phi(F)/|F|$. Let $\mathcal F=(F_n)$ be a sequence of finite subsets of $G$, and define 
		\[
		C(\mathcal F)=\sup_{E\Subset G}\limsup_{n\to\infty} \frac{|EF_n|}{|F_n|} \ \in [1,\infty].
		\]
		If  $C(\mathcal F)<\infty$, or $\mathfrak{h}>0$, or both, then we have \begin{equation}\label{eq:folnerness}\sup_{\phi\in\Phi}\limsup_{n\to\infty}\frac{\phi(F_n)}{|F_n|} = C(\mathcal F)\cdot \mathfrak{h}\end{equation}
		If  $C(\mathcal F)=\infty$ and also $\mathfrak{h}=0$, then no information is provided. 
	\end{theorem}
    \begin{remark}
        One can interpret $C(\mathcal F)$  as a ``F\o lnerness''  constant. In fact, it is easy to see that $C(\mathcal F)$  equals 1 if and only if $\mathcal F$ is F\o lner. 
    \end{remark}
We will be interested in the following natural choices of $\Phi$ (see \Cref{sec:preliminaries} for definitions).
    \begin{itemize}
        \item For a topological dynamical $G$-system $(X,T)$, we let $\Phi$ be the collection of functions $\phi\colon \mathscr{F}(G)\to[0,\infty)$ of the form $F\mapsto \log N(\mathcal U^F)$, where $\mathcal U$ is an open cover of $X$. In this case $\mathfrak h_{\phi}=h_{top}(T,\mathcal U)$ and $\mathfrak h=h_{top}(T)$. \item For a measure-preserving $G$-system $(X,\mathcal B,\mu,T)$, we let $\Phi$ be the collection of functions of the form $F\mapsto H_{\mu}(\xi^F)$, where $\xi$ is a finite measurable partition. Then $\mathfrak h_{\phi}=h_{\mu}(T,\xi)$ and $\mathfrak h=h_{\mu}(T)$. 
    \end{itemize}
It is well-known that for the computation of entropy (either topological or measure-theoretic), if we use a non-F\o lner sequence in the definition, then we may obtain something different from entropy. \Cref{thm:main-intro} shows that in the positive entropy case, we in fact obtain entropy scaled by the factor $C(\mathcal F)\in[1,\infty]$. In particular, every choice of a non-F\o lner sequence overestimates entropy, where the excess takes the form of a multiplicative factor  $C(\mathcal F)$. Thus this result can be seen as a version of the Ornstein-Weiss Lemma or Fekete's Lemma for non-F\o lner sequences, with the caveat that it concerns a family of sub-additive functions instead of individual functions. 
    
    Let us also mention that an elementary special case of \Cref{thm:main-intro} is the standard subgroup entropy formula for entropy: if $H\leqslant G$ is a subgroup and  $\mathcal F$ is a F\o lner sequence on $H$, then $C(\mathcal F)$ is equal to the index of $H$ in $G$.

	If $A=(t_i)_{i\in\N}$ is a sequence of elements in $G$ with no repetitions, then taking $F_n=\{t_1,\dots,t_n\}$ in \Cref{thm:main-intro} we recover  the known proportionality results for sequence entropy of Krug, Newton, and Goodman, and their  generalization to countable amenable groups  by Gao, Huang, and Zhu \cite{goodman_topological_1974,krug_sequence_1972,gao_sequence_2020}. If repetitions are allowed, however, then this is not the case anymore. Our main result (\Cref{thm:main}) is slightly more general than \Cref{thm:main-intro}, and also handles sequences with repetitions.
	
	Our proof relies on the existence of exact tilings of countable amenable groups \cite{downarowicz_tilings_2019}, and results about the infimum's rule from \cite{kolyada_shearers_2016}. These results ensure that the condition in the second paragraph of \Cref{thm:main-intro} is verified by the entropy functions we consider. Our results abstract those in \cite{carrasco-vargas_topological_2025}. We also prove a version of Goodwyn's theorem for sequence entropy in \Cref{sec:goodwyn}.
\section{Preliminaries}\label{sec:preliminaries}
	\subsection{Countable amenable groups}
	Here we review some facts about countably infinite amenable groups and functions on their finite subsets. The reader is referred to \cite[Chapters 4.5 and 4.7]{kerr_ergodic_2016} and \cite{kolyada_shearers_2016}.
	
	We fix in what follows a countably infinite amenable group $G$ with neutral element $1_G$. We write $E\Subset G$ when $E$ is a nonempty finite subset of $G$. The collection of nonempty finite subsets of $G$ will be denoted $\mathscr F({G})$.
	
	We employ the notations $EF=\{ef : e\in E, \ f\in F\}$ and  $Eg=\{eg : e\in E\}$, $E,F\Subset G$, $g\in G$. We say that $F$ is $(E,\epsilon)$-invariant when $|EF\triangle F|/|F|\leq \epsilon$, $\epsilon>0$. Here $\triangle$ denotes symmetric difference and $|\ \cdot \ |$ denotes cardinality. 
	
	Recall that a countable group is amenable if and only if for all $\epsilon>0$ and $E\Subset G$ there exists a finite subset of $G$ which is $(E,\epsilon)$-invariant. A sequence $(F_n)_{n\in\N}$ of finite subsets of $G$ is called a (left) F\o lner sequence when for every $\epsilon>0$ and $E\Subset G$ the set $F_n$ is eventually $(E,\epsilon)$-invariant.
	
	\begin{definition}\label{set-functions}
		Let $\phi\colon\mathscr F({G})\to[0,\infty)$ be a function. 
		\begin{enumerate}
			\item $\phi$ is monotone when $\phi(E)\leq \phi(F)$ for all $E\subset F$.
			\item $\phi$ is (right) translation-invariant if $\phi(Ef)=\phi(E)$ for all $E\Subset G$, $f\in G$.
			\item $\phi$ is subadditive if $\phi(E\cup F)\leq \phi(E)+\phi(F)$ for all $E,F\Subset G$. 
			\item $\phi(F)$ converges to a limit $\mathfrak{h}\in [0,\infty)$ as $F$ becomes more and more invariant if for all $\epsilon>0$ one can find $E\Subset G$ and $\delta>0$ such that 
			\[|\phi(F)-\mathfrak{h}|<\epsilon\]
			for any $F\Subset G$ which is $(E,\delta)$-invariant.
			\item $\phi$ converges to its infimum, or satisfies the infimum's rule, whenever the previous item holds and furthermore
			\[\inf_{F\Subset G}\phi(F)=\mathfrak{h}\]
		\end{enumerate}
	\end{definition}
    \begin{remark}
        If $\phi$ goes from $\mathscr{F}(G)$ to $[0,\infty)$, then its normalization  $F\to\phi(F)/|F|$ is also a function from $\mathscr{F}(G)$ to $[0,\infty)$.  \Cref{set-functions} can be applied to both of them.
    \end{remark}
	The next result is the analogue of Fekete's lemma on subadditive sequences for amenable groups, and is generally referred to as Ornstein-Weiss lemma as it originates in \cite{ornstein_entropy_1987}. This formulation is taken from \cite{kerr_ergodic_2016}. 
\begin{proposition}\label{fekete}
		If $\phi\colon\mathscr{F}(G)\to[0,\infty)$ is a translation-invariant and subadditive function then  
		\[F\mapsto \frac{1}{|F|}\phi(F)\]
		converges to a limit as $F$ becomes more and more invariant. 
	\end{proposition}
	See \cite[Theorem 4.38]{kerr_ergodic_2016} for a proof. The definition of $(E,\epsilon)$ invariant that we took is slightly different from the one considered in this reference (\cite[Definition 4.32]{kerr_ergodic_2016}), but we can interchange them thanks to \cite[Lemma 2.6]{downarowicz_tilings_2019}.

	\begin{remark}
		If $\phi$ is subadditive and translation-invariant then $\phi(E)\leq |E|\phi(\{1_G\})$ for all $E\Subset G$. More generally, we have $\phi(EF)\leq |F|\phi(E)$ for all $E,F\Subset G$. It suffices to write $EF=\bigcup_{f\in F}Ef$.
	\end{remark}
    
\subsection{Topological dynamical systems}
	By a topological dynamical system $(X,T)$ we mean that $X$ is a compact metrizable topological space and $T\colon G\times X\to X$, $(g,x)\mapsto T^g(x)$ is a left group action such that every map $x\mapsto T^g(x)$ is a homeomorphism.  
	
	Given a finite open cover $\mathcal U$ of $X$ we denote by $N(\mathcal U)$  the minimal cardinality of a subset of $\mathcal U$ which covers $X$. We write 
	\[\mathcal{U}^F=\bigvee_{g\in F}(T^{g})^{-1}(\mathcal U), \ \ F\Subset G.\]
	One easily verifies that for $E,F\Subset G$ we have $(\mathcal U^E)^F=\mathcal U^{EF}$. We also remark that the map 
	\[\mathscr{F}(G)\to[0,\infty), \ F\mapsto \log N(\mathcal U^{F})\]
	is subadditive and (right) translation-invariant as defined in \Cref{set-functions}.
	\begin{definition}
		The topological entropy $h_{top}(T,\mathcal U)$ of $T$  with respect to an open cover $\mathcal U$ is defined as  the limit of 
		\[\mathscr{F}(G)\to[0,\infty), \ F\mapsto \frac{1}{|F|}\log N(\mathcal U^{F})\]
		as $F$ becomes more and more invariant. The limit exists  by \Cref{fekete}. The topological entropy $h_{top}(T)$ of $T$ is defined as the supremum of $h_{top}(T,\mathcal U)$, with $\mathcal{U}$ ranging over all finite open covers of $X$. 
	\end{definition}
	For fixed  $\mathcal U$ the function $F\mapsto \log N(\mathcal U^F)/|F|$ converges to a limit $h_{top}(T,\mathcal U)$, but it is unknown if this limit equals the infimum of $\log N(\mathcal U^F)/|F|$ with $F$ ranging over finite subsets of $G$. However, the following result is proved in \cite[Theorem 6.8]{kolyada_shearers_2016}.
	\begin{proposition}\label{prop:infimum-rule-topological} The topological entropy $h_{top}(T)$ satisfies
		\[h_{top}(T)=\sup_{\mathcal U}\inf_{F\Subset G}\frac{1}{|F|}\log N(\mathcal{U}^F)\]
	\end{proposition}
	Next we define topological sequence entropy.
	\begin{definition}
		Let $A=(t_i)_{i\in\N}$ be a sequence of elements in $G$. The topological sequence entropy $h_{top}^A(T,\mathcal U)$ with respect to an open cover $\mathcal U$ is defined as 
		\[{h_{top}^{A}}(T,\mathcal U)=\limsup_{n\to\infty}\frac{1}{n}\log N(\mathcal U^{\{t_1,\dots,t_n\}})\]
		The topological sequence entropy $h_{top}^A(T)$ of $T$ with respect to the sequence $A$ is defined as the supremum of $h_{top}^{A}(T,\mathcal U)$, where $\mathcal U$ ranges over all finite open covers of $X$. 
	\end{definition}
	\subsection{Measure-preserving dynamical systems}
	By a measure-preserving dynamical system $(X,\mathcal X,\mu, T)$ we mean that $(X,\mathcal X,\mu)$ is a probability space and $T\colon G\times X\to X$, $(g,x)\mapsto T^g(x)$ is a left group action such that every map $x\mapsto T^g(x)$ is measurable and measure-preserving.
	
	Given a finite measurable partition $\xi$ of $X$ we denote
	\[H_{\mu}(\xi)=-\sum_{P\in\xi}\mu(P)\log \mu(P)\]
	We write 
	\[\xi^F=\bigvee_{g\in F}T^{g^{-1}}(\xi), \ \ F\Subset G.\]
	One easily verifies that for $E,F\Subset G$ we have $(\xi^E)^F=\xi^{EF}$. We remark that the map 
	\[\mathscr{F}(G)\to[0,\infty), \ F\mapsto H_{\mu}(\xi^{F})\]
	is subadditive and (right) translation-invariant as defined in \Cref{set-functions}.
	\begin{definition}
		The measure-theoretic entropy $h_{\mu}(T,\xi)$ of $T$ with respect to $\xi$ is defined as the limit of 
		\[\mathscr{F}(G)\to[0,\infty), \ F\mapsto \frac{1}{|F|}H_{\mu}(\xi^{F})\]
		as $F$ becomes more and more invariant. The limit exists  by \Cref{fekete}. The measure-theoretic entropy $h_{\mu}(T)$ of $T$ is defined as the supremum of $h_{\mu}(T,\xi)$, with $\xi$ ranging over all finite measurable partitions of $X$. 
	\end{definition}
	The following result is proved in \cite{kolyada_shearers_2016}, see also  \cite[Chapter 9.3]{kerr_ergodic_2016}.
	\begin{proposition}\label{prop:infimum-rule-measurable} Let $\xi$ be a finite measurable partition of $X$. The function $F\mapsto H_{\mu}(\xi^F)/|F|$ converges to its infimum, meaning that
		\[h_{\mu}(T,\xi)=\inf_{F\Subset G}\frac{1}{|F|}H_{\mu}(\xi^F)\]
	\end{proposition}
	Next we define measure-theoretic sequence entropy.
	\begin{definition}
		Let $A=(t_i)_{i\in\N}$ be a sequence of elements in $G$. The measure-theoretic sequence entropy $h_{\mu}^A(T,\xi)$ with respect to a finite measurable partition $\xi$ is defined as 
		\[{h_{\mu}^{A}}(T,\xi)=\limsup_{n\to\infty}\frac{1}{n}H_{\mu}(\xi^{\{t_1,\dots,t_n\}})\]
		The measure-theoretic sequence entropy $h_{\mu}^A(T)$ of $T$ with respect to the sequence $A$ is defined as the supremum of $h_{\mu}^{A}(T,\xi)$, where $\xi$ ranges over all finite measurable partitions of $X$. 
	\end{definition}
	\section{A result on subadditive functions}\label{sec:subadditive-functions}
	In this section we prove results about subadditive functions over finite subsets of $G$. These can be seen as abstract versions of classic results of Krug, Newton, and Goodman about entropy functions \cite{newton_sequence_1970,newton_sequence_1970a,krug_sequence_1972,goodman_topological_1974} (and more recently \cite{gao_sequence_2020}).  The main result of this section is \Cref{thm:main} below.\begin{definition}\label{constant-C-S} 
    Let $\mathcal F=(F_n)_{n\in\N}$ be a sequence of finite nonempty subsets of $G$, and let $\mathcal S=(s(n))_{n\in\N}$ be a sequence of positive real numbers.
		We define 
		\[C_E(\mathcal F,\mathcal S)=\limsup_{n\to\infty}\frac{|EF_n|}{s(n)}, \ \ E\Subset G.\]
        Next, we define 
\[C(\mathcal F,\mathcal S)=\sup_{E\Subset G}C_E(\mathcal F,\mathcal S) \ \ \in [0,\infty]
\]
The analogous quantities obtained by replacing $\limsup$ by $\liminf$ are denoted $\underline{C}_E(\mathcal F,\mathcal S)$ and $\underline{C}(\mathcal F,\mathcal S)$.
	\end{definition}
\begin{lemma}\label{upper-bound-versatile}
		Let $\phi\colon\mathscr{F}(G)\to[0,\infty)$ be a monotone translation-invariant subadditive function, and let  $\mathfrak{h}$ be the limit of $\phi(F)/|F|$ as $F$ becomes more and more invariant (which exists by  \Cref{fekete}). Let $\mathcal F$ and $\mathcal S$ be as in \Cref{constant-C-S}, and assume that $C(\mathcal F,\mathcal S)<\infty$, or $\mathfrak{h}> 0$, or both. Then
		\[\limsup_{n\to\infty}\frac{1}{s(n)}\phi(F_n)\leq C(\mathcal F,\mathcal S)\cdot \mathfrak{h}\]
		If $C(\mathcal F,\mathcal S)=0$ and $\mathfrak{h}=\infty$ then the right hand side is interpreted as $0$. If  $C(\mathcal F,\mathcal S)=\infty$ and also $\mathfrak{h}=0$ then no information is provided. 
		
		A similar result holds for the liminf of $\phi(F_n)/s(n)$ and $\underline{C}(\mathcal F,\mathcal S)$.
	\end{lemma}
	\begin{proof}
		We start with the claim about limsups. We will consider different cases, let us begin with the trivial ones. First observe that our claim is trivially true if one of $\{C(\mathcal F,\mathcal S), \mathfrak{h}\}$ is infinite and the other one is nonzero. Next, we consider the case in which $C(\mathcal F,\mathcal S)$ is zero (and $\mathfrak{h}$ could be infinite or not). Since $C(\mathcal F,\mathcal S)$ is zero, we have that $C_{\{1_G\}}(\mathcal F,\mathcal S)$ is also zero. By definition, this means that the limsup of $|F_n|/s(n)$ is zero as $n\to\infty$. Since $\phi(F_n)\leq |F_n|\phi(\{1_G\})$ and $\phi(\{1_G\})<\infty$, we have 
		\[\limsup_{n\to\infty}\frac{\phi(F_n)}{s(n)}\leq \limsup_{n\to\infty}\frac{|F_n|}{s(n)}\phi(\{1_G\})=0.\]
		Thus our claim is also verified in this case. Next, let us suppose that both $C(\mathcal F,\mathcal S)$ and $\mathfrak{h}$ are finite (this overlaps with the previous case if $C(\mathcal F,\mathcal S)$ is zero, but assuming that it is positive makes no difference in the next argument). This case is the interesting one. We prove that for all $\epsilon>0$ we have
		\begin{equation}\label{eq:upper-bound-versatile}
			\limsup_{n\to\infty}\frac{1}{s(n)}\phi(F_n)\leq C(\mathcal F,\mathcal S)\cdot (\mathfrak{h}+\epsilon)+\epsilon\cdot \phi(\{1_G\}).
		\end{equation}
		The claim in the statement follows from this inequality.
		
		From now on $\epsilon>0$ is fixed. Since $\phi(F)/|F|$ converges to a limit as $F$ becomes more and more invariant, we can find  $\delta>0$ and $W\Subset G$ such that every set $F\Subset G$ which is $(W,\delta)$-invariant satisfies
		\[\phi(F)\leq |F|(\mathfrak{h}+\epsilon).\]
		By \cite[Theorem 4.3]{downarowicz_tilings_2019} it is possible to choose a finite set  $\textbf{T}$ of $(W,\delta)$-invariant finite subsets of $G$ and sequences $(g_i)_{i\in\N}\in G^\N$ and $(T_i)_{i\in\N}\in \textbf{T}^\N$,  so that 
		\[G=\bigsqcup_{i\in\N} T_ig_i. \]
		Let $R=\bigcup_{U\in \textbf{T}}UU^{-1}$. We now choose an increasing sequence $(j(n))_{n\in\N}$ of natural numbers such that 
		\[	\limsup_{n\to\infty}\frac{1}{s(n)}\phi(F_n)=\lim_{n\to\infty}\frac{1}{s(j(n))}\phi(F_{j(n)}).\]
		Note that any subsequence of $(j(n))_{n\in\N}$ will also have this property, as a subsequence of a convergent sequence converges to the same limit. 
		Let us define
		\[M=\sup_{D\Subset G, \ 1_G\in D}\limsup_{n\to\infty}\frac{|DF_{j(n)}|}{s(j(n))}\]
		Observe that $M\leq C(\mathcal F,\mathcal S)$, and the second quantity is finite by hypothesis. Thus we can choose $Q\Subset G$ containing $1_G$ so that   
		\[
		\limsup_{n\to\infty} \frac{|QF_{j(n)}|}{s(j(n))}+\epsilon\geq M\]
		Choose a subsequence $(h(n))_{n\in\N}$ of  $(j(n))_{n\in\N}$ such that 
		\[\limsup_{n\to\infty} \frac{|QF_{j(n)}|}{s(j(n))}=\lim_{n\to\infty} \frac{|QF_{h(n)}|}{s(h(n))}.\]
		Observe that we have the inequalities
		\begin{equation}\label{eq:satanic-inequalities-limsups}
			M-\epsilon\leq \limsup_{n\to\infty} \frac{|QF_{j(n)}|}{s(j(n))} = \lim_{n\to\infty} \frac{|QF_{h(n)}|}{s(h(n))}\leq \limsup_{n\to\infty} \frac{|RQF_{h(n)}|}{s(h(n))}\leq \limsup_{n\to\infty} \frac{|RQF_{j(n)}|}{s(j(n))}\leq  M.
		\end{equation}
		The first inequality is by our choice of $Q$, the next one is by definition of $(h(n))_{n\in\N}$, the next one is because $Q\subset RQ$, the next inequality is because $h(n)$ is a subsequence of $j(n)$, and the last one is by definition of $M$. It will be convenient to let $E=RQ$.

		We now start with the computation of our  upper bound. We have
		\[
		\limsup_{n\to\infty}\frac{1}{s(n)}\phi(F_n)=\limsup_{n\to\infty}\frac{1}{s(j(n))}\phi(F_{j(n)})=\limsup_{n\to\infty}\frac{1}{s(h(n))}\phi(F_{h(n)})\]
		Since $1_G\in E$ and $\phi$ is monotone, we can bound the last term as
		\[
		\leq \limsup_{n\to\infty}\frac{1}{s(h(n))}\phi(EF_{h(n)})
		\]
		We now partition 
		$EF_{h(n)} = \bigsqcup_{i\in\N} EF_{h(n)}\cap T_ig_i$. Since this is actually a finite union and $\phi$ is subadditive we have  $\phi(EF_{h(n)})\leq \sum_{i\in\N}\phi(EF_{h(n)}\cap T_ig_i)$ (defining $\phi(\emptyset)=0$ to account for the empty intersections). The computation continues as  
		\[\leq \limsup_{n\to\infty}\sum_{i\in\N}\frac{1}{s(h(n))}\phi(EF_{h(n)}\cap T_ig_i)\]
		We now group these elements in a convenient manner. We say that $(n,i)$ is good when $EF_{h(n)}\cap T_ig_i$ is equal to $T_ig_i$,  and otherwise we say that $(n,i)$ is bad. Similarly, we say that $g\in EF_{h(n)}$ is good when it belongs to $EF_{h(n)}\cap T_ig_i$ and $(n,i)$ is good, and otherwise $g$ is bad.
		
		The intuitive idea is that the good elements in $EF_{h(n)}$ are the ones that live in a ``translated tile'' $T_ig_i$ which is fully covered by $EF_{h(n)}$; we call this good as this ensures that $g$ lives in a subset of $EF_{h(n)}$ which is sufficiently invariant. This should become clear in the rest of the proof. 
		
		We split the sum in the ``good part'' and the ``bad part'' as
		\[=\limsup_{n\to\infty}\left(\sum_{\substack{i\in\N\\ (n,i)\text{ good}}}\frac{1}{s(h(n))}\phi(EF_{h(n)}\cap T_ig_i)+\sum_{\substack{i\in\N\\ (n,i)\text{ bad}}}\frac{1}{s(h(n))}\phi(EF_{h(n)}\cap T_ig_i)\right)\]
		We will now bound the good part from above. Note that for good $(n,i)$ the set $EF_{h(n)}\cap T_ig_i=T_ig_i$ is $(W,\delta)$-invariant thanks to our choice of $\textbf{T}$. This ensures that we have the very desirable bound
		\[\phi(EF_{h(n)}\cap T_ig_i)\leq|EF_{h(n)}\cap T_ig_i| (\mathfrak{h}+\epsilon)\]
		Thus we can bound the sum of the good part as
		\begin{align*}
			\sum_{\substack{i\in\N\\ (n,i)\text{ good}}}\frac{1}{s(h(n))}\phi(EF_{h(n)}\cap T_ig_i)\leq &{\sum_{\substack{i\in\N\\ (n,i)\text{ good}}}\frac{1}{s(h(n))}|EF_{h(n)}\cap T_ig_i| (\mathfrak{h}+\epsilon)}\\
			=& \frac{|\{g\in EF_{h(n)} : g \ \text{good}\}|}{s(h(n))}(\mathfrak{h}+\epsilon)\\
			\leq&   \frac{|EF_{h(n)}|}{s(h(n))}(\mathfrak{h}+\epsilon)
		\end{align*}
		Thus the limsup of the good part can be bounded as
		\begin{align*}
			\limsup_{n\to\infty}\sum_{\substack{i\in\N\\ (n,i)\text{ good}}}\frac{1}{s(h(n))}\phi(EF_{h(n)}\cap T_ig_i)\leq& \limsup_{n\to\infty}\frac{|EF_{h(n)}|}{s(h(n))}(\mathfrak{h}+\epsilon)\\
			\leq & C_E(\mathcal F,\mathcal S)(\mathfrak{h}+\epsilon)\\
			\leq & C(\mathcal F,\mathcal S)(\mathfrak{h}+\epsilon)
		\end{align*}
		
		We now bound the bad part. For this we will exploit our choice of $E$ and the subsequences we took. The idea is that when $(n,i)$ is bad then we can not control how invariant is  $EF_{h(n)}\cap T_ig_i$.  Instead, our choice of $E$ ensures that most elements in $EF_{h(n)}$ are in the good part. We remark that here we strongly rely on the hypothesis $C(\mathcal F,\mathcal S)<\infty$.
		
		By subadditivity and translation-invariance we have $\phi(F)\leq |F|\phi(\{1_G\})$, and therefore
		\[\sum_{\substack{i\in\N\\ (n,i)  \ \text{bad}}}\frac{1}{s(h(n))}\phi(EF_{h(n)}\cap T_ig_i)\leq \sum_{\substack{i\in\N\\ (n,i)  \ \text{bad}}}\frac{|EF_{h(n)}\cap T_ig_i|}{s(h(n))}\phi(\{1_G\})=\frac{|\{g \in EF_{h(n)} : g \ \text{bad}\}|}{s(h(n))}\phi(\{1_G\})
		\]
		Recall that $E=RQ$, that $R=\bigcup_{U\in \textbf{T}} UU^{-1}$, and that every $T_i$ belongs to $\textbf{T}$. We can now explain  our choice of $R$. Take an arbitrary element $g=tg_i$ belonging to $T_ig_i$ for some $i\in\N$. Then an arbitrary element $t_ig_i$ in $T_ig_i$ can be written as $t_ig_i=t_it^{-1}tg_i=t_it^{-1}g\in T_iT^{-1}_ig\subset Rg$. Therefore for any $g\in T_ig_i$ we have $T_ig_i\subset Rg$.
		
		We claim that every element from $EF_{h(n)}$ which belongs to $QF_{h(n)}\subset EF_{h(n)}$ is good. Indeed, let $g\in QF_{h(n)}$ be arbitrary. By the observation in the previous paragraph the translated tile $T_ig_i$ to which $g$ belongs must be contained in $Rg$, which is contained in $RQF_{h(n)}=EF_{h(n)}$. That is, $T_ig_i$ is fully covered by $EF_{h(n)}$, which means that $g$ is good by definition.

		Note that $1_G\in R$, which implies that $QF_{h(n)}$ is contained in $RQF_{h(n)}$, which then implies $|RQF_{h(n)}\smallsetminus QF_{h(n)}|=|RQF_{h(n)}|-|QF_{h(n)}|$. Therefore we have
		\[
		\frac{|\{g \in EF_{h(n)}: g \ \text{bad}\}|}{s(h(n))}\leq \frac{|RQF_{h(n)}\smallsetminus QF_{h(n)}|}{s(h(n))}=\frac{|RQF_{h(n)}|}{s(h(n))}-\frac{|QF_{h(n)}|}{s(h(n))}
		\]
		Taking limsup we obtain
		\begin{align*}
			{\limsup_{n\to\infty} \frac{|\{g \in EF_{h(n)} : g \ \text{bad}\}|}{s(h(n))}}\leq& {\limsup_{n\to\infty}\left(\frac{|RQF_{h(n)}|}{s(h(n))}-\frac{|QF_{h(n)}|}{s(h(n))}\right)}\\
			=&{\limsup_{n\to\infty} \frac{|RQF_{h(n)}|}{s(h(n))} - \liminf_{n\to\infty} \frac{|QF_{h(n)}|}{s(h(n))}}
		\end{align*}
		Our choice of $(h(n))_{n\in\N}$ ensures that the liminf is actually a limit, so this value is equal to 
		\[={\limsup_{n\to\infty} \frac{|RQF_{h(n)}|}{s(h(n))} - \lim_{n\to\infty} \frac{|QF_{h(n)}|}{s(h(n))}}\]
		Thanks to \Cref{eq:satanic-inequalities-limsups}, we can bound this quantity as 
		\[
		\leq \epsilon.
		\]
		Therefore the limsup of the bad part is at most $\epsilon\phi(\{1_G\})$. We already proved  that the limsup of the good part is at most $C(\mathcal F,\mathcal S)(\mathfrak{h}+\epsilon)$, so combining these two we obtain \Cref{eq:upper-bound-versatile}.  
		
		We now prove the statement regarding liminfs. Let us assume $\underline{C}(\mathcal F,\mathcal S)<\infty$ and $\mathfrak{h}<\infty$, as the other cases can be handled with the same argument from before. It suffices to prove that for all $\epsilon>0$ we have
		\begin{equation}\label{eq:liminf-upper-bound-versatile}\liminf_{n\to\infty}\frac{\phi(F_n)}{s(n)}\leq (\underline{C}(\mathcal F,\mathcal S)+\epsilon)(\mathfrak{h}+\epsilon)+\epsilon\phi(\{1_G\}).\end{equation}
		The proof is very similar to the proof of \Cref{eq:upper-bound-versatile}, so we briefly explain the modifications required. 
		
		Liminfs and limsups behave differently by taking subsequences, so we must modify the order of our choices to adapt the previous argument. After fixing $\epsilon$ one chooses $(W,\delta)$, $\textbf{T}$, and $R$ in the same manner as before. Then we choose $Q\Subset G$ containing $1_G$ so that $\underline{C}_Q(\mathcal F,\mathcal S)+\epsilon\geq \underline{C}(\mathcal F,\mathcal S)$, and we let $E=RQ$. Then we choose an increasing sequence of natural numbers $(h(n))_{n\in\N}$ so that 
		\[\liminf_{n\to\infty}\frac{|EF_{n}|}{s(n)}=\lim_{n\to\infty}\frac{|EF_{h(n)}|}{s(h(n))}\]
		Observe that
		\begin{equation}\label{eq:satanic-versatile-upper-bound-liminfs}
			\underline{C}(\mathcal F,\mathcal S)-\epsilon\leq \liminf_{n\to\infty}\frac{|QF_n|}{s(n)}\leq \liminf_{n\to\infty}\frac{|QF_{h(n)}|}{s(h(n))}\leq \liminf_{n\to\infty}\frac{|EF_{h(n)}|}{s(h(n))}=\liminf_{n\to\infty}\frac{|EF_{n}|}{s(n)}\leq \underline{C}(\mathcal F,\mathcal S)
		\end{equation}
		The first inequality is by our choice of $Q$, the next one is because $(h(n))_{n\in\N}$ is a subsequence of $(n)_{n\in\N}$, the next one is because $Q\subset E$, the next one is by our choice of $(h(n))_{n\in\N}$, and the last one is by definition of $\underline{C}(\mathcal F,\mathcal S)$.
		
		Now the proof of \Cref{eq:liminf-upper-bound-versatile} is as follows. Since the liminfs are non-decreasing by taking subsequences, we have
		\[\liminf_{n\to\infty}
		\frac{1}{s(n)}\phi(F_n)\leq \liminf_{n\to\infty}
		\frac{1}{s(h(n))}\phi(F_{h(n)})\]
		We define ``good'' and ``bad'' in the same way as before. The term at the right is split in good and bad using sub-additivity and we get
		\[\leq\liminf_{n\to\infty}\left(\sum_{\substack{i\in\N\\ (n,i)\text{ good}}}\frac{1}{s(h(n))}\phi(EF_{h(n)}\cap T_ig_i)+\sum_{\substack{i\in\N\\ (n,i)\text{ bad}}}\frac{1}{s(h(n))}\phi(EF_{h(n)}\cap T_ig_i)\right)\]
		Recall that for any pair of sequences $(a_n)_{n\in\N}$ and $(b_n)_{n\in\N}$ we have $ \liminf_{n\to\infty}(a_n+b_n)\leq \liminf_{n\to\infty}a_n+\limsup_{n\to\infty}b_n$. Thus we can continue the computation as
		\[\leq \liminf_{n\to\infty}\sum_{\substack{i\in\N\\ (n,i)\text{ good}}}\frac{1}{s(h(n))}\phi(EF_{h(n)}\cap T_ig_i)+\limsup_{n\to\infty}\sum_{\substack{i\in\N\\ (n,i)\text{ bad}}}\frac{1}{s(h(n))}\phi(EF_{h(n)}\cap T_ig_i).
		\]
		The same argument from before shows that the liminf of the good part is at most 
		\[\underline{C}_E(\mathcal F,\mathcal S)(\mathfrak{h}+\epsilon)\leq \underline{C}(\mathcal F,\mathcal S)(\mathfrak{h}+\epsilon).\]
		Furthermore, the same argument from before shows that the limsup of the bad part is at most 
		\[\limsup_{n\to\infty}\frac{|\{g\in EF_{h(n)} :  g\text{ bad}\}|}{s(h(n))}\phi(\{1_G\})\]
		With the same argument as before, our choice of $E=RQ$ allows us to bound
		\[\limsup_{n\to\infty}\frac{|\{g\in EF_{h(n)} :  g\text{ bad}\}|}{s(h(n))}\leq {\limsup_{n\to\infty} \frac{|RQF_{h(n)}|}{s(h(n))} - \liminf_{n\to\infty} \frac{|QF_{h(n)}|}{s(h(n))}}\]
		Now, thanks to our choice of $(h(n))_{n\in\N}$, the first limsup is in fact a limit, so the last quantity is equal to
		\[{=\liminf_{n\to\infty} \frac{|RQF_{h(n)}|}{s(h(n))} - \liminf_{n\to\infty} \frac{|QF_{h(n)}|}{s(h(n))}}\]
		This value is at most $\epsilon$ by \Cref{eq:satanic-versatile-upper-bound-liminfs}, and this finishes the argument.
	\end{proof}
    \begin{theorem}\label{thm:main}
		Let $\Phi$ be a family of monotone subadditive and translation-invariant functions from $\mathscr{F}(G)$ to $[0,\infty)$. Suppose that $\Phi$ has the following closure property: for every $\phi\in\Phi$ and $E\Subset G$ the function $F\to\phi(EF)$ 
		belongs to $\Phi$. The function defined in this way is subadditive and translation-invariant. 
        
        For each $\phi\in \Phi$ let $\mathfrak{h}_\phi$ be the limit of $\phi(F)/|F|$ as $F$ becomes more and more invariant, which exists by \Cref{fekete}, and let 
		\[\mathfrak{h}=\sup_{\phi\in\Phi}\mathfrak{h}_{\phi}.\]
		We assume that the following weak form of the infimum's rule holds:
		\[\mathfrak{h}=\sup_{\phi\in\Phi}\inf_{F\Subset G}\frac{\phi(F)}{|F|}.\]
		Let $\mathcal F$ and $\mathcal S$ be as in \Cref{constant-C-S}, and assume that  $C(\mathcal F,\mathcal S)<\infty$, or $\mathfrak{h}>0$, or both. Then
		\[\sup_{\phi\in\Phi}\limsup_{n\to\infty}\frac{\phi(F_n)}{s(n)} = C(\mathcal F,\mathcal S)\cdot \mathfrak{h}\]
		This quantity is interpreted as $0$ if $C(\mathcal F,\mathcal S)=0$ and $\mathfrak{h}=\infty$. If  $C(\mathcal F,\mathcal S)=\infty$ and also $\mathfrak{h}=0$ then no information is provided. 
		
		The analogous statement holds if we replace $\limsup$ by $\liminf$ and instead consider $\underline{C}(\mathcal F,\mathcal S)$.
	\end{theorem}

	\begin{proof}
		We only prove the claim regarding $\limsup$ as the other case is very similar. It will be convenient to re-write the statement as
		\begin{equation}\label{eq:familiar-theorem}
			\sup_{\phi\in\Phi}\limsup_{n\to\infty}\frac{1}{s(n)}\phi(F_n) =\begin{cases}
				C(\mathcal F,\mathcal S)\cdot \mathfrak{h} &C(\mathcal F,\mathcal S)<\infty, \mathfrak{h}<\infty\\
				\infty \ &C(\mathcal F,\mathcal S)>0, \mathfrak{h}=\infty\\
				\infty &C(\mathcal F,\mathcal S)=\infty, \mathfrak{h}>0\\
				0 &C(\mathcal F,\mathcal S)=0, \mathfrak{h}=\infty
		\end{cases}\end{equation}
		We start by proving the inequality $\leq$ in \Cref{eq:familiar-theorem}. Let us first consider the case in which $C(\mathcal F,\mathcal S)<\infty$ and $\mathfrak{h}<\infty$. Let $\phi\in\Phi$, and observe that $\mathfrak{h}_{\phi}<\infty$ because $\mathfrak{h}<\infty$. Then \Cref{upper-bound-versatile} shows that
		\[\limsup_{n\to\infty}\frac{\phi(F_n)}{s(n)}\leq C(\mathcal F,\mathcal S)\cdot \mathfrak{h}_{\phi}.\]
		Taking supremum over $\phi$ we obtain the desired inequality. Next, let us consider the case in which $C(\mathcal F,\mathcal S)=0$ and $\mathfrak{h}=\infty$. For $\phi\in\Phi$,  \Cref{upper-bound-versatile} asserts that
		\[\limsup_{n\to\infty}\frac{\phi(F_n)}{s(n)}\leq 0.\]
		Taking supremum over $\phi$ we obtain the inequality $\leq$ in \Cref{eq:familiar-theorem}. The remaining two cases hold for trivial reasons. 
		
		We now prove the inequality $\geq$ in \Cref{eq:familiar-theorem}. If $\mathfrak{h}=0$ then $C(\mathcal F,\mathcal S)<\infty$ by hypothesis, and then the inequality
		\[
		\sup_{\phi\in\Phi}\limsup_{n\to\infty}\frac{1}{s(n)}\phi(F_n)\geq C(\mathcal F,\mathcal S)\cdot \mathfrak{h}=0
		\]
		holds for trivial reasons. Suppose now that $\mathfrak{h}>0$. Let us temporarily say that $\phi\in\Phi$ is good when 
		\[
		\inf_{D\Subset G}\frac{\phi(D)}{|D|}>0
		\]
		Observe that the assumption that $\mathfrak{h}>0$ and the weak form of the infimum's rule that we assume imply that the collection of good elements in $\Phi$ is nonempty. Also note that for every $\phi\in\Phi$ and $F\Subset G$ we have 
		\[
		\phi(F)\geq |F|\cdot\inf_{D\Subset G}\frac{\phi(D)}{|D|}
		\]
		Fix $E\Subset G$. Applying the closure property of $\Phi$, the previous observation, the definition of $C_E(\mathcal F,\mathcal S)$, and the weak form of the infimum's rule, we have
		\begin{align*}
			\sup_{\phi\in\Phi}\limsup_{n\to\infty}\frac{1}{s(n)}\phi(F_n) &\geq \sup_{\phi\in\Phi}\limsup_{n\to\infty}\frac{1}{s(n)}\phi(EF_n) \\
			&\geq \sup_{\phi\text{ good}}\limsup_{n\to\infty}\frac{1}{s(n)}\phi(EF_n)\\
			&\geq \sup_{\phi\text{ good}}\limsup_{n\to\infty}\left(\frac{1}{s(n)}\cdot |EF_n|\cdot\inf_{D\Subset G}\frac{\phi(D)}{|D|}\right)\\
			&=\sup_{\phi\text{ good}}\left(C_E(\mathcal F,\mathcal S)\cdot\inf_{D\Subset G}\frac{\phi(D)}{|D|}\right)\\
			&=C_E(\mathcal F,\mathcal S)\cdot \sup_{\phi \text{ good}}\inf_{D\Subset G}\frac{\phi(D)}{|D|}\\
			&= C_E(\mathcal F,\mathcal S)\cdot \mathfrak{h}.
		\end{align*}
		Since this is true for arbitrary $E$ we can take supremum over $E\Subset G$ and obtain
		\[\sup_{\phi\in\Phi}\limsup_{n\to\infty}\frac{1}{s(n)}\phi(F_n)\geq C(\mathcal F,\mathcal S)\cdot \mathfrak{h}.\]
		This finishes the argument.
	\end{proof}
    
	\section{The theorems of Krug, Newton and Goodman for amenable groups}\label{sec:sequence-entropy}
	In this section we  derive from \Cref{thm:main} the mentioned results about the proportionality principle for topological and measure-theoretic sequence entropy. 
\begin{definition}\label{constant-K}
		Let $A=(t_n)_{n\in\N}$ be a sequence of elements in $G$. We define
		\[K_E(A)=\limsup_{n\to\infty}\frac{|E\{t_1,\dots,t_n\}|}{n}, \ E\Subset G.\]
		Then we define $K(A)$ as the supremum of $\{K_E(A) : E\Subset G\}$. 
	\end{definition}
	This corresponds to taking $\mathcal F$ as $(\{t_1,\dots,t_n\})_{n\in\N}$, and taking $\mathcal S=(s(n))_{n\in\N}$ as $s(n)=n$ in \Cref{constant-C-S}. This is the classical constant introduced in Newton's work \cite{newton_sequence_1970}, and subsequently considered in \cite{krug_sequence_1972,goodman_topological_1974,gao_sequence_2020}. We will also be interested in the following constant, which we partially defined in the statement of \Cref{thm:main-intro}.
\begin{definition}\label{constant-C}
		Let $\mathcal F=(F_n)_{n\in\N}$ be a  sequence of finite subsets of $G$. We define \[C_E(\mathcal F)=\limsup_{n\to\infty}\frac{|EF_n|}{|F_n|}, \ E\Subset G.\]
		Then we define $C(\mathcal F)$ as the supremum of $\{C_E(\mathcal F) : E\Subset G\}$. The analogous quantity with $\liminf$ is denoted  $\underline{C}(\mathcal F)$.
	\end{definition}
	This corresponds to choosing $\mathcal S=(s(n))_{n\in\N}$ as $s(n)=|F_n|$ in \Cref{constant-C-S}. We find this a natural choice as it can be interpreted as a ``F\o lnerness constant'' of the sequence.
	\begin{proposition}\label{folner-sequence-iff-c-equals-one}
		$\mathcal F$ is a (left) F\o lner sequence if and only if $C(\mathcal F)=1$.
	\end{proposition}
	\begin{proof}
		One easily verifies that $\mathcal F=(F_n)_{n\in\N}$ is a (left) F\o lner sequence if and only if for all $E\Subset G$ with $1_G\in E$ we have $|EF_n\triangle F_n|/|F_n|\to 0$ as $n\to\infty$, and that $C(\mathcal F)$ equals 
		the supremum of $\{C_E(\mathcal F) : E\Subset G, \ 1_G\in E\}$. Now fix $E\Subset G$ with $1_G\in E$. Observe that $|EF_n\triangle F_n|=|EF_n\smallsetminus F_n|=|EF_n|-|F_n|$, and therefore
		$|EF_n|/|F_n|=1+|EF_n\triangle F_n|/|F_n|$. Taking the limsup of this expression as $n\to\infty$ we see that $C_E(\mathcal F)=1$ if and only if $|E F_n\triangle F_n|/|F_n|\to 0$ as $n\to\infty$. 
	\end{proof}
	\subsection{Topological sequence entropy}\label{subsec:topological-sequence-entropy}
	Let $(X,T)$ be a topological dynamical system (in the sense of \Cref{sec:preliminaries}). 
    
    \begin{proposition}\label{prop:topological-entropy-functions}
    Let $\Phi$ be the collection of functions $\mathscr{F}(G)\to[0,\infty)$ of the form $F\to\log N(\mathcal U^F)$, where $\mathcal U$ is a finite open cover of $X$. Then $\Phi$ satisfies the hypotheses in \Cref{thm:main}.
    \end{proposition}
    \begin{proof}
    The relation $(\mathcal U^E)^F=\mathcal U^{EF}$ shows that $\Phi$ has the closure property required in \Cref{thm:main}. Furthermore, the weak form of the infimum's rule is verified by \Cref{prop:infimum-rule-topological}.
    \end{proof}
\begin{theorem}\label{thm:sequence-entropy-rescaling}
		Let $A=(t_i)_{i\in\N}$ be a sequence of elements in $G$. If either $K(A)<\infty$ or $h_{top}(T)>0$ then we have
		\[h^A_{top}(T)=K(A)\cdot h_{top}(T).\]
		The term at the right is interpreted as zero if $K(A)=0$ and $h_{top}(T)=\infty$. If $K(A)=\infty$ and $h_{top}(T)=0$ then no information is provided.
	\end{theorem}
	\begin{proof}
Apply \Cref{thm:main} with $\Phi$ as in \Cref{prop:topological-entropy-functions}, $\mathcal F=(F_n)_{n\in\N}$ given by $F_n=\{t_1,\dots,t_n\}$, and $\mathcal S=(s(n))_{n\in\N}$ given by $s(n)=n$.
\end{proof}

	We also prove the following result, which corresponds to the computation of topological entropy through a sequence of finite sets not required to be F\o lner, and which we believe may have independent interest.
	\begin{theorem}\label{topological-entropy-non-folner-sequence}
		Let $\mathcal F=(F_n)_{n\in\N}$ be a sequence of finite subsets of $G$, and let $C(\mathcal F)$ be the quantity defined in \Cref{constant-C}. If either $C(\mathcal F)<\infty$ or $h_{top}(T)\ne 0$ then 
		\[\sup_{\mathcal U}\limsup_{n\to\infty}\frac{\log N(\mathcal U^{F_n})}{|F_n|}=C(\mathcal F)\cdot h_{top}(T).\]
		Here $\mathcal U$ ranges over all finite open covers of $X$. The analogue statement holds with $\liminf$ and the constant $\underline{C}(\mathcal F)$.
	\end{theorem}
	\begin{proof}
Apply \Cref{thm:main} with $\Phi$ as in \Cref{prop:topological-entropy-functions} and $\mathcal S=(s(n))_{n\in\N}$ given by $s(n)=|F_n|$.
\end{proof}
	The previous result shows that we can interchange limsup by liminf exactly when $\underline{C}(\mathcal F)=C(\mathcal F)$, and that otherwise the gap between them is determined by the gap between $\underline{C}(\mathcal F)$ and $C(\mathcal F)$ ($0<h_{top}(T)<\infty$). Since $C(\mathcal F)=1$ exactly when $\mathcal F$ is F\o lner (\Cref{folner-sequence-iff-c-equals-one}), we have the following  corollary of \Cref{topological-entropy-non-folner-sequence}. 
\begin{corollary}
		Let  $\mathcal F=(F_n)_{n\in\N}$ be a sequence of finite subsets of $G$ which is not a F\o lner sequence, let $C(\mathcal F)$ be the quantity defined in \Cref{constant-C}, and suppose that $0<h_{top}(T)<\infty$. Then  
		\[\sup_{\mathcal U}\limsup_{n\to\infty}\frac{\log N(\mathcal U^{F_n})}{|F_n|}>h_{top}(T).\]
		Here $\mathcal U$ ranges over all finite open covers of $X$.
	\end{corollary}
	We note that this was  proved for $G=\Z$ in \cite[\S 4.1]{carrasco-vargas_topological_2025}.
    \subsection{Measure-theoretic sequence entropy}
	
	Let $(X,\mathcal X,\mu, T)$ be a measure-preserving dynamical system in the sense of \Cref{sec:preliminaries}. 
        \begin{proposition}\label{prop:measurable-entropy-functions}
    Let $\Phi$ be the collection of functions $\mathscr{F}(G)\to[0,\infty)$ of the form $F\to H_{\mu}(\xi^F)$, where $\xi$ is a finite measurable partition of $X$. Then $\Phi$ satisfies the hypotheses in \Cref{thm:main}.\end{proposition}
    \begin{proof}
        The relation $(\xi ^E)^F=\xi^{EF}$ shows that $\Phi$ has the closure property required in \Cref{thm:main}. Furthermore, the condition about the infimum's rule is verified by \Cref{prop:infimum-rule-measurable}.
    \end{proof}\begin{theorem}\label{thm:measure-sequence-entropy-rescaling}
		Let $A=(t_i)_{i\in\N}$ be a sequence of elements in $G$. If either $K(A)<\infty$ or $h_{\mu}(T)>0$ then we have
		\[h^A_{\mu}(T)=K(A)\cdot h_{\mu}(T).\]
		The term at the right is interpreted as zero if $K(A)=0$ and $h_{\mu}(T)=\infty$. If $K(A)=\infty$ and $h_{\mu}(T)=0$ then no information is provided.
	\end{theorem}
	
	\begin{proof}
Apply \Cref{thm:main} with $\Phi$ as in \Cref{prop:measurable-entropy-functions}, $\mathcal F=(F_n)_{n\in\N}$ given by $F_n=\{t_1,\dots,t_n\}$, and $\mathcal S=(s(n))_{n\in\N}$ given by $s(n)=n$.
\end{proof}

	We also have the following result, which corresponds to computing measure-theoretic entropy through a sequence which is not required to be a F\o lner sequence.
	\begin{theorem}\label{measure-entropy-non-folner-sequence}
		Let $\mathcal F=(F_n)_{n\in\N}$ be a sequence of finite subsets of $G$, and let $C(\mathcal F)$ be the quantity defined in \Cref{constant-C}. If either $C(\mathcal F)<\infty$ or $h_{\mu}(T)\ne 0$ then 
		\[\sup_{\xi}\limsup_{n\to\infty}\frac{H_{\mu}(\xi^{F_n})}{|F_n|}=C(\mathcal F)\cdot h_{\mu}(T).\]
		Here $\xi$ ranges over all finite measurable partitions $X$. The analogue statement holds with $\liminf$ and the constant $\underline{C}(\mathcal F)$.
	\end{theorem}
	
	\begin{proof}
Apply \Cref{thm:main} with $\Phi$ as in \Cref{prop:measurable-entropy-functions} and $\mathcal S=(s(n))_{n\in\N}$ given by $s(n)=|F_n|$.
\end{proof}
	The previous result shows that we can interchange limsup by liminf exactly when $\underline{C}(\mathcal F)=C(\mathcal F)$, and that otherwise the gap between them is given by the gap between $\underline{C}(\mathcal F)$ and $C(\mathcal F)$. Since $C(\mathcal F)=1$ exactly when $\mathcal F$ is F\o lner (\Cref{folner-sequence-iff-c-equals-one}), we have the following amusing corollary of \Cref{measure-entropy-non-folner-sequence}. \begin{corollary}
		Let  $\mathcal F=(F_n)_{n\in\N}$ be a sequence of finite subsets of $G$ which is not a F\o lner sequence, and suppose that $0<h_{\mu}(T)<\infty$. Then  
		\[\sup_{\xi}\limsup_{n\to\infty}\frac{H_{\mu}(\xi^{F_n})}{|F_n|}>h_{\mu}(T).\]
		Here $\xi$ ranges over all finite measurable partitions of $X$.
\end{corollary}

\section{Goodwyn's theorem for sequence entropy}\label{sec:goodwyn}
While the variational principle is not always satisfied by sequence entropy \cite{goodman_topological_1974}, the measure-theoretic sequence entropy is always bounded by the topological one. In this section we prove the following result.
\begin{theorem}\label{thm:goodwyn-sequence}
    Let $(X,T)$ be a topological dynamical system and let $A=(t_i)_{i\in\N}$ be a sequence of elements in $G$. Then
    \[\sup_{\mu}h^A_{\mu}(T)\leq h^A_{top}(T),\]
    where the supremum is taken over all $T$-invariant Borel probability measures on $X$. Furthermore, if either $K(A)<\infty$ or $h_{top}(T)\ne 0$, then we have an equality.
\end{theorem}


Let $\xi$ be a finite measurable partition of $X$ and let $\mu$ be a $T$-invariant Borel probability measure on $X$. We denote by $N(\xi,\delta)$ the minimal number of atoms in $\xi$ needed to cover a subset of $X$ with measure at least $1-\delta$, where $0<\delta<1$. The first step in our proof is the next elementary proposition. 
\begin{proposition}
    Let $0<\delta<1/e$. For every $F\Subset G$ we have 
    \[
H_{\mu}(\xi^F)\leq \log N(\xi^F,\delta)+\delta |F|\log(|\xi|)+(\log 2+2/e)
    \]
\end{proposition}
\begin{proof}

Let $\gamma$ be a collection of atoms in $\xi^F$ whose union has measure at least $1-\delta$ and minimal cardinality $N(\xi^F,\delta)$.  Let $P$ be the union of the atoms in $\gamma$, and let $Q=X\smallsetminus P$. Let $\alpha$ be the partition of $X$ defined by adding to $\gamma$ the atom  $Q$, and let $\beta$ be the partition of $X$ consisting of the atom $P$, and all atoms in $\xi^F$ that are disjoint from $P$. It follows that $\xi^F=\alpha\vee\beta$ and consequently $H_{\mu}(\xi^F)\leq H_{\mu}(\alpha)+H_{\mu}(\beta)$. 

Since $|\alpha|\leq N(\xi^F,\delta)+1\leq 2N(\xi^F,\delta)$, we have
\[H_{\mu}(\alpha)\leq\log|\alpha|\leq \log (2N(\xi^F,\delta))\leq \log N(\xi^F,\delta)+\log 2.\]
On the other hand, Proposition 3.3 in \cite{carrasco-vargas_topological_2025} shows that  $H_{\mu}(\beta)\leq 2/e+\delta \log(|\beta|)$, which can be bounded by $2/e+\delta|F|\log(|\xi|)$ using the trivial bound $|\beta|\leq|\xi^F|\leq |\xi|^{|F|}$. The sum of both upper bounds equals the right hand side of the inequality in the statement  (alternatively, one can reach a similar conclusion by conditioning to the sigma-algebra generated by $\{P,Q\}$, and using \cite[fact 1.6.1]{downarowicz_entropy_2011}). 
\end{proof}
Next we upper bound $N(\xi^F,\delta)$ using covering numbers associated to Hamming pseudometrics. Endow the finite set $\Omega(\xi,F)$ of functions from $F$ to $\xi$ with the Hamming distance \[d_{\Omega(\xi,F)}^H(u,v)=\frac{1}{|F|}|\{i\in F : u(i)\ne v(i)\}|.\] 
The cardinality of an $\epsilon$-ball in $d_{\Omega(\xi,F)}^H$ does not depend on its center, so we simply denote it $|B_{\Omega(\xi,F)}(\epsilon)|$. 

For $x$ and $y$ in $X$, we define $d^H_{\xi,F}(x,y)$ as the Hamming distance $d^H_{\Omega(\xi,F)}$ of the finite codes $(\xi(T^g(x)))_{g\in F}$ and $(\xi(T^g(y)))_{g\in F}$ (where $\xi(x)$ denotes the atom containing $x$). That is 
\[d^H_{\xi,F}(x,y)=\frac{1}{|F|}|\{g\in F : T^g(x)\text{ and }T^g(y)\text{ are in different atoms of }\xi\}|\]
We denote by $\spa^H_{\xi}(F,\epsilon,\delta)$ the minimal number of $\epsilon$-balls in $d^H_{\xi,F}$ needed to cover measure $\mu$ at least $1-\delta$.  
\begin{proposition}\label{prop:N-and-hamming}
    For $\epsilon,\delta\in(0,1)$ and $F\Subset G$ we have 
\begin{equation}
    N(\xi^F,\delta)\leq \spa_\xi^H(F,\epsilon,\delta) \cdot |B_{\Omega(\xi,F)}(\epsilon)|
\end{equation}
\end{proposition}
\begin{proof}It suffices to observe that an $\epsilon$-ball in  $d^H_{\xi,F}$ equals the union of at most $|B_{\Omega(\xi,F)}(\epsilon)|$ atoms of $\xi^F$. 
\end{proof}
When we take logarithms  and divide by $|F|$ (or a bigger number), the contribution of $|B_{\Omega(\xi,F)}(\epsilon)|$ becomes negligible as $|F|\to\infty$ and $\epsilon\to 0$. 
\begin{proposition}[Equation 1.3 in \cite{katok_lyapunov_1980}]\label{prop:katok-bound}
    For $\epsilon<(|\xi|-1)/|\xi|$ we have 
\begin{equation}
\lim_{|F|\to\infty}\frac{1}{|F|}\log(|B_{\Omega(\xi,F)}(\epsilon)|)=\epsilon\log(|\xi|-1)-\epsilon\log\epsilon-(1-\epsilon)\log(1-\epsilon)\end{equation}
\end{proposition}
We will use Bowen covering numbers to upper bound $\spa^H_\xi(F,\epsilon,\delta)$. Let $d$ be a compatible metric for $X$, and for $F\Subset G$ define the (Bowen) metric $d_F^B$ by \[d^B_F(x,y)=\max\{d(T^g(x),T^g(y)) : g\in F\}.\]

Let $\spa^B_d(F,\epsilon,\delta)$ be the minimal number of $\epsilon$-balls in $d^B_F$ needed to cover a subset of $X$ with measure $\mu$ at least $1-\delta$, $\delta\in(0,1)$.

For $P\in \xi$, define $\partial_\delta P$ as $\{x\in P : d(x,y)<\delta\text{ for some } y\in X\smallsetminus P\}$. Furthermore, let $\partial_\delta \xi=\bigcup_{P\in\xi}\partial_\delta P$. Denoting by $\partial \xi$ the union of the topological boundaries of atoms in $\xi$, we have
\[
\partial \xi=\bigcap_{\delta>0}\partial_\delta\xi
\]
It follows that whenever $\partial\xi$ has zero measure, the measure of $\partial_\delta \xi$ tends to zero as $\delta\to 0$. 
Here $\partial\xi$ denotes the union of the topological boundaries of its atoms.   
\begin{proposition}[Proposition 2 in \cite{katok_slow_1997}]\label{prop:katok-thouvenot-bound}
    Suppose that $\mu(\partial\xi)=0$. Let $0<\delta<1/2$, and let $0<\lambda<1$ be small enough so that $\mu(\partial_\lambda \xi)<\delta^2$. Then
    \[
\spa^H_{\xi}(F,\delta,2\delta)\leq\spa^B_{d}(F,\lambda/2,\delta)
    \]
\end{proposition}
We remark that \cite{katok_slow_1997} actually assumes $\mu$ non-atomic, while we do not. Inspecting the proof in \cite{katok_slow_1997} it is not difficult to see that the inequality remains valid even if $\mu$ has atoms. Furthermore, note that we took a different definition of $\partial_\lambda$  (following \cite{katok_lyapunov_1980}).

Let $\spa^B_d(F,\epsilon)$ be the minimal number of $\epsilon$-balls in $d^B_F$ needed to cover $X$. Putting the previous results together, we obtain an upper bound for $H_{\mu}(\xi^F)$ in terms of $\spa^B_d(F,\epsilon)$ (and error terms that will vanish). This upper bound plays the same role as the one below equation 6.4 in \cite{kolyada_shearers_2016}, but here the measure does not need to be ergodic.

\begin{proposition}\label{prop:horrible-inequality-hamming}
Suppose that $\mu(\partial \xi)=0$. Let $0<\delta<1/(2e)$ and let $0<\lambda<1$ be small enough so that $\mu(\partial_\lambda \xi)<\delta^2$. Then for all $F\Subset G$ we have
    \[H_{\mu}(\xi^F)\leq \log (\spa^B_d(F,\lambda/2))+\log(|B_{\Omega(\xi,F)}(\delta)|)+2\delta|F|\log(|\xi|)+ (\log 2 + 2/e)\]
\end{proposition}
\begin{proof}
Applying \Cref{prop:N-and-hamming} and then \Cref{prop:katok-thouvenot-bound} we obtain
\begin{align*}
    H_{\mu}(\xi^F)&\leq \log N(\xi^F,2\delta)+2\delta|F|\log(|\xi|)+(\log2+2/e)\\
    &\leq \log \spa^H_{\xi}(F,\delta,2\delta)+\log(|B_{\Omega(\xi,F)}(\delta)|)+2\delta|F|\log(|\xi|)+(\log2+2/e)\\
    &\leq \log \spa^B_{d}(F,\lambda/2,\delta)+\log(|B_{\Omega(\xi,F)}(\delta)|)+2\delta|F|\log(|\xi|)+(\log2+2/e)
\end{align*}
Then observe that $\spa^B_{d}(F,\lambda/2,\delta)\leq \spa^B_{d}(F,\lambda/2)$. 
\end{proof}
The standard argument shows that topological sequence entropy, which we defined with open covers, can also be computed as
\[
h_{top}^A(T)=\lim_{\epsilon\to 0}\limsup_{n\to\infty}\frac{1}{n}\log \spa^B_d(\{t_1,\dots,t_n\},\epsilon)
\]

\begin{proof}[Proof of \Cref{thm:goodwyn-sequence}]
    Let $\mu$ be a $T$-invariant Borel probability measure on $X$. Let $\xi$ be a finite Borel measurable partition of $X$ such that $\mu(\partial \xi)=0$. Let $0<\delta<1/(2e)$, and let $\lambda(\delta)$  satisfy $\mu(\partial_{\lambda(\delta)}\xi)<
    \delta^2$, $0<\lambda(\delta)<1$. Applying \Cref{prop:horrible-inequality-hamming} with $F=\{t_1,\dots,t_n\}$ and dividing by $n$, which is at least the cardinality of $\{t_1,\dots,t_n\}$, we obtain
\begin{equation}
    \limsup_{n\to\infty}\frac{1}{n}H_{\mu}(\xi^{\{t_1,\dots,t_n\}})\leq h^A_{top}(T)+f(\delta)
\end{equation}
where $f$ is the error term and \Cref{prop:katok-bound} ensures $\lim_{\delta\to0}f(\delta)=0$. Taking $\delta\to 0$ we obtain 
\[h^A_{\mu}(T,\xi)\leq h^A_{top}(T).\]
To finish the argument, we note that $h^A_{\mu}(T)$ equals the supremum of $h^A_{\mu}(T,\xi)$ where $\xi$ ranges over finite Borel measurable partitions of $X$ with $\mu(\partial\xi)=0$, as these partitions can be taken with arbitrarily small diameter \cite[Lemma 8.5]{walters_introduction_2000}. This is the case thanks to \cite[Lemma 2]{kushnirenko_metric_1967}, which is also valid in the case of amenable group actions via the same proof. Thus we have proved inequality $\leq$ in \Cref{thm:goodwyn-sequence}. The fact that we have an equality whenever $K(A)<\infty$ or $h_{top}(T)>0$ follows from  \Cref{thm:sequence-entropy-rescaling},  \Cref{thm:measure-sequence-entropy-rescaling}, and the standard variational principle. 
\end{proof}

\begin{remark}\label{rem:previous-proofs}
Previous proofs \cite{goodman_topological_1974,eberlein_topological_1976,huang_combinatorial_2009} of \Cref{thm:goodwyn-sequence} for actions of a single transformation invoke Lemma 4 in Kushnirenko's paper  \cite{kushnirenko_metric_1967}, which was later found to be false as stated \cite{lemanczyk_sequence_1985}. To our understanding, the first work where this gap is filled is \cite{gao_sequence_2020}. Our proof follows a different path from those in \cite{goodman_topological_1974,eberlein_topological_1976,huang_combinatorial_2009,gao_sequence_2020}, and instead generalizes the argument in \cite[\S 3]{carrasco-vargas_topological_2025}.
\end{remark}
\section*{Acknowledgements}
Part of this work was done during the Simons Semester ``Invariant Structures in Groups, Spaces, and Algebras'', at IMPAN, Poland, and the author acknowledges financial support from the Simons Foundation grant (award no. SFI-MPS-T-Institutes-00010825), and State Treasury funds as part of a task commissioned by the Minister of Science and Higher Education under the project “Organization of the Simons Semesters at the Banach Center - New Energies in 2026-2028” (agreement no. MNiSW/2025/DAP/491). Furthermore, this work was supported by the ERC grant DynAMiCs (101167561) of the European Research Council.
\bibliography{references}
\end{document}